\documentclass{amsart}

\RequirePackage{hyperref}

\usepackage{amsmath,amsfonts,tikz,tikz-3dplot,capt-of,hyperref,enumitem,mathtools}
\usepackage{graphicx}
\newtheorem{lem}{Lemma}
\newtheorem{thm}[lem]{Theorem}

\newtheorem{prop}[lem]{Proposition}
\newtheorem{defn}[lem]{Definition}
\newtheorem{rmk}[lem]{Remark}

\newtheorem{ex}[lem]{Example}

\numberwithin{lem}{section}

\begin{document}

\title[Local Gromov-Witten invariants of some snc Surfaces]{Local Gromov-Witten Invariants of Some Simple Normal Crossing Surfaces}

\author{Sheldon Katz}
\address{Department of Mathematics, University of Illinois, Urbana IL 61801}
\email{katzs@illinois.edu}
\author{Sungwoo Nam}
\address{Department of Mathematics, University of Illinois, Urbana IL 61801}
\email{sungwoo2@illinois.edu}

\begin{abstract}
Inspired by M-theory and superconformal field theory, we extend the notions of local Gromov--Witten invariants from the case of del Pezzo surfaces to shrinkable surfaces, a class of reducible surfaces with simple normal crossings satisfying certain positivity conditions.
\end{abstract}

\maketitle

\section{Introduction}

In physics, 5-dimensional superconformal field theories (5d SCFTs) can be realized starting from a noncompact Calabi--Yau threefold $X$ admitting an analytic contraction $\pi:X\to Y$ to a small neighborhood of a canonical threefold singularity $(Y,p)$ (not necessarily isolated) \cite{Mosei,XieYau}.   The rank of the gauge group of the SCFT is the number of compact irreducible 2-dimensional components of $S=f^{-1}(p)$.  In the case of rank~1, if $S$ is smooth then in a neighborhood of $S$ the threefold $X$ is analytically isomorphic to a neighborhood of the zero section of a local del Pezzo surface.  

In \cite{JKKV}, the above notion of a contractible threefold was generalized to the notion of a \emph{shrinkable threefold}.  Shrinkable threefolds also realize 5d SCFTs, the rough idea being that volumes of surfaces can be made to approach zero by taking an appropriate limit in K\"ahler moduli space, effectively implementing a contraction.  

Since 5d SCFTs have BPS invariants, we anticipate that local Gromov-Witten invariants can be associated to shrinkable threefolds of arbitrary rank, generalizing local Gromov-Witten invariants of del Pezzo surfaces in the rank~1 case.  We show that this intuition borrowed from physics is correct. 
 Under certain hypotheses, we can define local Gromov-Witten invariants (Theorem~\ref{prop:localSiebert} and (\ref{eqn:localGW})) in terms of $S$ alone.  Our main results are Theorems~\ref{thm:main} and~\ref{thm:generalcase}.  Theorem~\ref{thm:main} shows that under certain hypotheses related to shrinkability, these local Gromov-Witten invariants are the contributions of $S$ to the Gromov-Witten invariants of any Calabi-Yau threefold $X$ containing $S$.  
 After extending the notion of local Gromov-Witten invariants to a more general context, Theorem~\ref{thm:generalcase} shows that under certain more general hypotheses related to shrinkability, these local Gromov-Witten invariants are the contributions of a scheme set-theoretically supported on $S$ to the Gromov-Witten invariants of any Calabi-Yau threefold $X$ containing $S$.

The paper \cite{JKKV} also contained a conjectural criterion for shrinkability which is more amenable to geometric analysis.  In fact, in most of \cite{JKKV}, this conjecture was assumed and was effectively used as an alternate definition of shrinkability.  Using that analysis, a classification of rank two 5d SCFTs was given up to a notion of physical equivalence, resulting in 64 physical equivalence classes.  In this paper, we only say a few words about physical equivalence, directing the interested reader to \cite{JKKV}.

We will abuse physics terminology slightly by defining a mathematical notion of shrinkability which is based on the physics conjecture.  We acknowledge that our definition only conjecturally agrees with the original physics definition, which we will say just a few words about below.

Let $X$ be a noncompact Calabi-Yau threefold which is an analytic neighborhood of a finite union $S=\cup S_i$ of projective surfaces.  We consider divisor classes $J=\sum a_i[S_i]$ on $X$, with integral $a_i > 0$.\footnote{In \cite{JKKV}, the $a_i$ were real numbers.  However, if a $J$ exists with real $a_i$ satisfying the conditions of Definition~\ref{def:shrinkable}, then since the intersection numbers $S_i\cdot C$ and $S_i\cdot S_j\cdot S_k$ are integers, we can replace the $a_i$ by rational numbers while maintaining these conditions.  Then we can clear denominators and assume that the $a_i$ are integers.}

\begin{defn}\label{defn:presh}
    $X$ is \emph{shrinkable} if there exists a $J$ with the following properties:
\begin{enumerate}[label=\emph{\textup{(}\roman*}\textup{)}]
    \item $J\cdot C\le 0$ for all curves $C\subset S$
    \item $J^2S_i\ge0$ for all $i$
    \item There exists at least one $i$ with $J^2S_i>0$.
\end{enumerate}
\label{def:shrinkable}
\end{defn}

\begin{rmk}\label{rmk:shdef}
Conditions (ii) and (iii) in Definition~\ref{defn:presh} will not be used again in this paper.  We have nevertheless included these conditions in the definition of a shrinkable surface to match the definition in \cite{JKKV}.  In cases where $J^2S_i=0$ for all $i$, the 5d quantum field theory resulting from M-theory compactification is related to a 6d theory compactified on a circle and is not an SCFT.  We give an example of such an $S$ immediately following the proof of Theorem~\ref{thm:generalcase}.
\label{rmk:extraneous}    
\end{rmk}

By abuse of notation, we follow \cite{JKKV} and say that a surface $S=\cup S_i$ is \emph{shrinkable} if $S$ can be embedded in a shrinkable threefold $X$.  We also refer to the number of irreducible components $S_i$ as the rank of $S$.

For example, if $X$ is the total space of the canonical bundle of a del Pezzo surface $S$, then we can use $J=a[S]$ for any $a>0$ to show that $X$ and hence $S$ is shrinkable.  M-theory compactification gives a rank one 5-dimensional gauge theory, and we arrive at a 5d SCFT by letting $a$ approach 0.  

However, the Hirzebruch surface $\mathbb{F}_2$ is a rank~1 shrinkable surface even though it is not a del Pezzo surface, as $X$ can be taken to be local $\mathbb{F}_2$, the total space of $\omega_{F_2}$, and the conditions of Definition~\ref{def:shrinkable} are then immediately checked.  This example may seem to contradict the previous paragraph.  However, this theory is physically equivalent to the $\mathbb{F}_0$ SCFT and gives nothing new.

\smallskip
Now assume that $S=\cup S_i$ has simple normal crossings (snc) singularities.  Since $\mathcal{O}_X(S_i)\vert_{S_i}\simeq\omega_{S_i}$ by adjunction, the intersection numbers in (\ref{def:shrinkable}) can be computed in terms of the geometry of $S$ alone, independent of the choice of Calabi-Yau $X$ containing $S$.  For example, if $C\subset S_i$ is a curve, then letting $C_{ij}=S_i\cap S_j$ we have
\begin{equation}
    \left(S_j\cdot C\right)_X = \left\{
    \begin{array}{cl}
       \left(K_{S_i}\cdot C\right)_{S_i}  & j=i \\
         \left(C_{ij}\cdot C\right)_{S_i}&j\ne i 
    \end{array}\right.
\label{eqn:computebyadjunction}
\end{equation}
and then $J\cdot C$ can be computed by linearity.  Since the right-hand side of (\ref{eqn:computebyadjunction}) makes sense without reference to $X$, we can define the expressions $J\cdot C$ using (\ref{eqn:computebyadjunction}) for the surface $S=\cup S_i$ without assuming the existence of an embedding of $S$ into a Calabi-Yau $X$.  Similarly, $J^2S_i$ can be formally defined as an intersection on $S_i$ without assuming the existence of a Calabi-Yau embedding.  

If $S=\cup S_i$ is snc and shrinkable, then there are some simple consequences. By the existence of $X$ we have
\begin{equation}
(C_{ij})^2_{S_i}+(C_{ij})^2_{S_j}=2g_{ij}-2,
\label{eqn:CYcondition}
\end{equation}
where $g_{ij}$ is the genus of $C_{ij}$.  We call (\ref{eqn:CYcondition}) the \emph{Calabi-Yau condition}.  Next, let $C\subset S_i$ be an irreducible curve not equal to any $C_{ij}$.   Then (1) in Definition~\ref{def:shrinkable} implies that $S_i\cdot C\le0$.  If in addition $C$ is a smooth rational curve, then $S_i\cdot C\le0$ is equivalent to $(C^2)_{S_i}\ge -2$.  In particular, $(-2)$ curves played an important role in the original definition of shrinkability in \cite{JKKV}, with noncompact ruled surfaces $N_k$ being glued to $S$ with a fiber of $N_k$ glued to $C$ (still contained inside $X$).  It was further conjectured in 
\cite{JKKV} that these hypotheses implied the existence of a contraction $X\to Y$, where $Y$ has canonical singularities, $S$ maps to a point $p$, and $\cup N_k$ contracts along fibers to a union of noncompact curves containing $p$.

\begin{defn}
An snc surface $S=\cup S_i$ is called \emph{pre-shrinkable} if it satisfies the Calabi-Yau condition (\ref{eqn:CYcondition}) and the conditions of Definition~\ref{def:shrinkable}, interpreted via (\ref{eqn:computebyadjunction}) and the subsequent discussion.
\end{defn}
 
 In their classification of rank 2 SCFTs, the authors of \cite{JKKV} implicitly assumed that pre-shrinkable surfaces are shrinkable, i.e.\ that $S$ can be embedded in a noncompact Calabi-Yau threefold. In this paper, we refer to the question of whether they are actually shrinkable (i.e.\ whether they can be embedded in a Calabi-Yau $X$) as the \emph{embeddability question}.  

If $J$ in Definition~\ref{def:shrinkable} can be taken to be $J=\sum [S_i]$ (so that $\omega_S^\vee$ is nef), then we can define $N^g_\beta\in\mathbb{Q}$ for $g\ge0$ and $\beta\in H_2(S,\mathbb{Z})$ satisfying a mild condition (Theorem~\ref{thm:main}), even if $S$ is only pre-shrinkable.  If $S$ is shrinkable and $S\subset X$, then these $N^g_\beta$ can be interpreted as local Gromov-Witten invariants, the contribution of stable maps to $S$ to the ordinary Gromov-Witten invariants of $X$.

The special cases of those $\beta$ not covered by Theorem~\ref{thm:main} are related to the physical notion of decoupled sectors used in the formulation of physical equivalence in \cite{JKKV}.

For general $J$ and shrinkable $S\subset X$, we define infinitesimal thickenings $S\subset S_J\subset X$ of $S$ in $X$ and describe local Gromov-Witten invariants in terms of the moduli stack of stable maps to $S_J$.  Theorem~\ref{thm:generalcase} shows that under certain more hypotheses, these local Gromov-Witten invariants are the contributions of $S_J$ to the Gromov-Witten invariants of any Calabi-Yau threefold $X$ containing $S_J$.

In \cite{Nam}, the second author has answered the embeddability question in the affirmative for 62 of the 64 physical equivalence classes in rank 2, completing the geometric construction of the resulting rank~2 SCFTs described in \cite{JKKV}.\footnote{This list does not include theories with $O7^+$-planes which are not known to admit a description using shrinkable surfaces.}.  In these cases, we have therefore rigorously formulated the local Gromov-Witten invariants of $S$ in terms of stable maps to infinitesimal thickenings of $S$. We give two examples of such shrinkable surfaces in Example~\ref{ex:preshrink}.  They are easily checked to be pre-shrinkable, and the embeddability question was answered affirmatively for these surfaces in \cite{Nam}.

\smallskip\noindent
{\bf Acknowledgements.}  We would like to thank Hee-Cheol Kim and Davesh Maulik for helpful conversations. The work of the authors was partially supported by NSF grants DMS 18-02242 and DMS 22-01203.

\section{Definition of Local Gromov--Witten Invariants}\label{sec:localgw}

In this section, we consider local Gromov--Witten invariants in a general setting.  Before setting up the problem, we pause to review Siebert's formula for the virtual fundamental class \cite{Siebert}.

Let $X$ be a smooth projective variety, and consider the universal stable map

\begin{equation}\label{eqn:univmap}
    \begin{array}{ccc}
    \overline{\mathcal{C}}_g(X,\beta)&\stackrel{f}{\rightarrow}&X\\
    \phantom{\pi}\downarrow\pi&&\\
    \overline{\mathcal{M}}_{g}(X,\beta)&&
    \end{array}
\end{equation}
Then the forgetful map $\rho:\overline{\mathcal{M}}_{g}(X,\beta)\to{\mathfrak{M}}_{g}$ supports a relative perfect obstruction theory $(R\pi_*f^*T_X)^\vee\to L^\bullet_\rho$, where ${\mathfrak{M}}_{g}$ is the Artin stack of prestable curves of genus $g$ \cite{BF}.  

On any scheme $Y$, we have a map
\begin{equation}
    \mathrm{ind}:\operatorname{Perf}(Y)\to K(Y)
\end{equation}
sending a perfect complex up to quasi-isomorphism to its K-theory class.  In particular $R\pi_*f^*T_X$ has a global resolution by a two-term complex of vector bundles $F_0\to F_1$ and we get a K-theory class
$\mathrm{ind}(R\pi_*f^*T_X) =[F_0]-[F_1]\in K(\overline{\mathcal{M}}_{g}(X,\beta))$.

Let $D=c_1(X)\cdot\beta+(\dim X -3)(1-g)$ be the virtual dimension of $\overline{\mathcal{M}}_{g}(X,\beta)$.  Then Siebert's formula is
\begin{equation}\label{eq:siebert}
    [\overline{\mathcal{M}}_{g}(X,\beta)]^{\mathrm{vir}}=\left(
    c\left(\mathrm{ind}(R\pi_*f^*T_X)\right)^{-1}\cap
    c_F\left(\overline{\mathcal{M}}_{g}(X,\beta)/{\mathfrak{M}}_{g}\right)\right)_D.
\end{equation}
In (\ref{eq:siebert}), $c_F\left(\overline{\mathcal{M}}_{g}(X,\beta)/{\mathfrak{M}}_{g}\right)$ is Fulton's relative canonical class for $\overline{\mathcal{M}}_{g}(X,\beta)$ relative to the Artin stack ${\mathfrak{M}}_{g}$ of prestable curves, and the subscript of $D$ means as usual that we only keep the $D$-dimensional contribution.

Now let $S=\cup S_i$ be any surface which is locally a hypersurface in a smooth threefold (for example, an snc surface). Let $\beta\in H_2(S,\mathbb{Z})$, and let $T_S^\bullet=R\underline{Hom}_{\mathcal{O}_S}(\Omega^1_S,\mathcal{O}_S)$ be the derived dual of $\Omega^1_S$.
\begin{lem}
The object $R\pi_*Lf^*T_S^\bullet\in D^b(\overline{\mathcal{M}}_g(S,\beta))$ is perfect.
\label{lem:perfect}
\end{lem}
In this lemma, $\pi$ and $f$ have the same meaning as they did in (\ref{eqn:univmap}), with $X$ replaced by $S$.

\begin{proof}
 We claim that we have a short exact sequence 
\begin{equation}
 0\rightarrow F_1\rightarrow F_0\rightarrow \Omega_S\rightarrow 0,
 \label{eqn:resomega}
\end{equation}
where $F_0$ and $F_1$ are locally free sheaves on $S$.  We have a surjection $\varphi:F_0\to \Omega_S$ with $F_0$ locally free since $\Omega_S$ is coherent on the projective variety $S$.  We put $F_1=\ker\varphi$ and show that $F_1$ is locally free, which we need only check on stalks at points $s\in S$.  The snc hypothesis implies that a neighborhood $U$ of $s\in S$ is a hypersurface in a smooth threefold $Y$.  Let $\mathcal{J}$ be the ideal sheaf of $U$ in $Y$.  Then the short exact sequence of $\mathcal{O}_{S,s}$-modules
\[
0\rightarrow \left(\mathcal{J}/\mathcal{J}^2\right)_s\rightarrow \Omega_{X,s}/\mathcal{J}\Omega_{X,s}\rightarrow \Omega_{S,s}\rightarrow 0
\]
shows that $\Omega_{S,s}$ has homological dimension~1, and then standard results of commutative algebra show that $(F_1)_s$ is a free $\mathcal{O}_{S,s}$-module.  

Since $T_S^\bullet$ is obtained by dualizing (\ref{eqn:resomega}), we deduce an exact triangle
\[
T_S^\bullet \rightarrow F_0^\vee\rightarrow F_1^\vee \stackrel{+}{\rightarrow},
\]
leading to the triangle
\begin{equation}
R\pi_*Lf^*T_S^\bullet \rightarrow R\pi_*f^*F_0^\vee\rightarrow R\pi_*f^*F_1^\vee \stackrel{+}{\rightarrow}.    
\label{eq:compobs}
\end{equation}
Since $F_0^\vee$ and $F_1^\vee$ are locally free, it follows that the second and third objects in the triangle (\ref{eq:compobs}) are perfect of amplitude $[0,1]$ \cite{BehrendGW}. It then follows that $R\pi_*Lf^*T_S^\bullet$ (which describes the relative obstruction theory for $\overline{\mathcal{M}}_{g}(S,\beta)\to{\mathfrak{M}}_{g})$ is also perfect, although not necessarily of amplitude $[0,1]$.   
\end{proof}

Now let $S\subset X$ be an snc hypersurface in a Calabi-Yau threefold $X$, with $i$ the inclusion.  Let $\beta\in H_2(S,\mathbb{Z})$.    
Let $T_\pi^\bullet$ be the derived dual of $\Omega_\pi$.  Suppose that stable maps to $S$ do not deform off $S$.  More precisely:

\begin{prop}\label{prop:localSiebert}
Suppose that $\overline{\mathcal{M}}_{g}(S,\beta)$ is a union of connected components of 
$\overline{\mathcal{M}}_{g}(X,i_*\beta)$ (as stacks).
Then the contribution of $\overline{\mathcal{M}}_{g}(S,\beta)$ to $[\overline{\mathcal{M}}_{g}(X,i_*\beta)]^{\mathrm{vir}}$ is
\begin{equation}
\left(c\left(-\mathrm{ind}\left(R\pi_*Lf^*T_S^\bullet\right)-\mathrm{ind}\left(R\pi_*f^*\omega_S\right)+\mathrm{ind}\left(R\pi_*T_\pi^\bullet\right)\right)\cap
    c_F\left(\overline{\mathcal{M}}_{g}(S,\beta)\right)\right)_0.
\label{eqn:localGWprop}
 \end{equation}
\end{prop}
Denoting this contribution by $[\overline{\mathcal{M}}_{g}(S,\beta)]^\mathrm{vir}$, we can define the local GW invariant as
\begin{equation}
N^g_\beta=\deg [\overline{\mathcal{M}}_{g}(S,\beta)]^\mathrm{vir},
\label{eqn:localGW}
\end{equation}
and $N^g_\beta$ is clearly independent of the choice of $X$ containing $S$ which satisfies the hypothesis of Proposition~\ref{prop:localSiebert}.  

\smallskip
We can define invariants $N^g_\beta$ for any pre-shrinkable snc surface $S$ using (\ref{eqn:localGWprop}) and (\ref{eqn:localGW}) without assuming the existence of any $X$.  However, without further hypotheses, these invariants should not be viewed as local Gromov-Witten invariants or as being associated with a perfect obstruction theory.

In Proposition~\ref{prop:components}, we will improve Proposition~\ref{prop:localSiebert} by giving criteria on $S$ and $\beta$ which ensures that the hypothesis of Proposition~\ref{prop:localSiebert} holds for any $X$ containing $S$.  That justifies calling  these $N^g_\beta$ local Gromov-Witten invariants of $S$, without assuming the existence of $X$, since they are contributions to the Gromov-Witten invariants of $X$ whenever $X$ exists.  In the majority of cases where the embeddability question has been answered in the affirmative~\cite{Nam}, the $N^g_\beta$ are indeed local Gromov-Witten invariants in the usual sense.

\begin{proof}
Let $\mathcal{I}$ be the ideal sheaf of $S$ in $X$.  After dualizing
\[
0\rightarrow \mathcal{I}/\mathcal{I}^2 \rightarrow \Omega^1_X|_S\rightarrow \Omega^1_S\rightarrow 0,
\]
pulling back by $f$ and pushing forward by $\pi$ we get
\begin{equation}\label{eq:indrel}
\mathrm{ind}\left(R\pi_*f^*T_X \right)=\mathrm{ind}\left(R\pi_*Lf^*T_S^\bullet\right)+\mathrm{ind}\left(R\pi_*f^*\omega_S\right),
\end{equation}

We also need the following formula from \cite{Siebert} for $c_F\left(\overline{\mathcal{M}}_{g}(S,\beta)/{\mathfrak{M}}_{g}\right)$.  
\begin{equation}
  c_F\left(\overline{\mathcal{M}}_{g}(S,\beta)/{\mathfrak{M}}_{g}\right)=  
  c\left(\mathrm{ind}\left(R\pi_*T_\pi^\bullet\right)\right)\cap c_F\left(\overline{\mathcal{M}}_{g}(S,\beta)\right).
\label{eqn:relabsFulton}
\end{equation}
The proposition follows immediately after substituting (\ref{eq:indrel}) and (\ref{eqn:relabsFulton}) into
(\ref{eq:siebert}).
\end{proof}

We next show that if $H^0(S,f^*\omega_S)=0$ for all $f\in\overline{\mathcal{M}}_g(S,\beta)$, the hypothesis of Proposition~\ref{prop:localSiebert} holds, so that local Gromov-Witten invariants can be defined. This condition is in turn implied by $J\cdot\beta<0$ with $J=\sum[S_i]$, a special case of the first condition for shrinkability in Definition~\ref{def:shrinkable}. 
 Here we have used $\mathcal{O}_X(J)\vert_S\simeq\omega_S$.
The condition on $H^0(f^*\omega_S)$ holds for all $\beta\in H_2(S,\mathbb{Z})$ if $\omega_S^\vee$ is ample.

\begin{prop}\label{prop:components}
Assume that $H^0(f^*\omega_S)=0$ for all $f\in\overline{\mathcal{M}}_g(S,\beta)$. 

Then $\overline{\mathcal{M}}_{g}(S,\beta)$ is a union of connected components of 
$\overline{\mathcal{M}}_{g}(X,i_*\beta)$ (as stacks).
\end{prop}

\begin{proof}  
The proof proceeds by applying the relative version of results of \cite{BF} to $S$ and $X$, and then comparing.

Let $T$ be a $\mathbb{C}$-scheme and 
consider any family of genus $g$ stable maps
\begin{equation}\label{eq:family}
    \begin{array}{ccc}
    \mathcal{C}&\stackrel{f}{\rightarrow}&X\\
    \phantom{\pi}\downarrow\pi&&\\
    T&&
    \end{array},
\end{equation}
with classifying map $g:T\to \overline{\mathcal{M}}_{g}(X,i_*\beta)$.  Put $h=\rho\circ g:T\to{\mathfrak{M}}_{g}$.

 Let $t_0\in T$ be a closed point and let $T\subset T'$. Let the ideals of $t_0$ in $T'$ and $T\subset T'$ be denoted by $\mathfrak{m}$ and $\mathcal{J}$ respectively.  Assume that $\mathfrak{m}\cdot \mathcal{J}=0$, which implies that $\mathcal{J}$ can be identified with a skyscraper sheaf on the reduced point $t_0$.
 
 Let $h':T'\to {\mathfrak{M}}_{g}$ be an extension of $h$ to $T'$ and put $C_{t_0}=\pi^{-1}(t_0)$.  Then there is an obstruction class 
\begin{equation}\label{eq:obst}
    \alpha\in H^1\left(C_{t_0},f^*T_X\right)\otimes \mathcal{J}
\end{equation}
whose vanishing is necessary and sufficient for the family of stable maps (\ref{eq:family}) to extend to a family over stable maps over $T'$
\begin{equation}
    \begin{array}{ccc}
    \mathcal{C}'&\stackrel{f'}{\rightarrow}&X\\
    \phantom{\pi'}\downarrow\pi'&&\\
    T'&&
    \end{array},
\end{equation}
whose classifying map $g':T'\to \overline{\mathcal{M}}_{g}(X,i_*\beta)$ satisfies $\rho\circ g' =h'$.
Furthermore, if $\alpha=0$, then the set of extensions is naturally a principal homogeneous space for $H^0\left(C_{t_0},f^*T_X\right)\otimes \mathcal{J}$.  The assertions in this paragraph are relative versions of the corresponding assertions in \cite{BF}.  The proofs of these assertions can be adapted to the relative situation in a straightforward manner. We will provide a little more detail about the relative obstruction theory presently in the more general setting of target varieties $S$ which are projective but not necessarily smooth.

By projectivity, we have a moduli stack of stable maps to $S$ and describe a relative obstruction theory 
$(R\pi_*Lf^*TS^\bullet)^\vee$
for $\overline{\mathcal{M}}_{g}(S,\beta)\to{\mathfrak{M}}_{g}$.  Here $T_S^\bullet$ is the derived dual $R\underline{\mathrm{Hom}}(L^\bullet_S,\mathcal{O}_S)\in D^{\mathrm{b}}(S)$.\footnote{For hypersurfaces $S\subset X$, we have $L^\bullet_S\simeq\Omega^1_S$.  Nevertheless, it will be convenient for generalizations to continue to use $L^\bullet_S$.}  We again explicitly adapt the results written in \cite{BF} to this situation.

Consider the universal stable map
\begin{equation}\label{eq:universal}
    \begin{array}{ccc}
    \overline{\mathcal{C}}_g(S,\beta)&\stackrel{f}{\rightarrow}&S\\
    \phantom{\pi}\downarrow\pi&&\\
    \overline{\mathcal{M}}_{g}(S,\beta)&&\\
    \end{array},
\end{equation}
which gives rise to the maps
\begin{equation}\label{eq:maps}
    Lf^*L^\bullet_S\rightarrow L^\bullet_{\overline{\mathcal{C}}_g(S,\beta)}\rightarrow L^\bullet_{\overline{\mathcal{C}}_g(S,\beta)/{\mathcal{C}}_g},
\end{equation}
where ${\mathcal{C}}_g$ is the universal curve over ${\mathfrak{M}}_{g}$.
We also have the fiber diagram
\begin{equation}\label{eq:fiber}
\begin{array}{ccc}
\overline{\mathcal{C}}_g(S,\beta)&\stackrel{\pi}{\rightarrow}&\overline{\mathcal{M}}_{g}(S,\beta)\\
\downarrow&&\downarrow\\
{\mathcal{C}}_g  &\rightarrow&{\mathfrak{M}}_{g}  
  \end{array},
\end{equation}
which induces an isomorphism $\pi^*(L^\bullet_{\overline{\mathcal{M}}_{g}(S,\beta)/{\mathfrak{M}}_{g}})\to L^\bullet_{{\mathcal{C}}_g(S,\beta)/{\mathcal{C}}_g}$.

Combining this isomorphism with (\ref{eq:maps}), we get a map
\begin{equation}
    Lf^*L^\bullet_S\rightarrow  \pi^*L^\bullet_{\overline{\mathcal{M}}_{g}(S,\beta)/{\mathfrak{M}}_{g}},
\end{equation}
which induces a map 
\begin{equation}
    Lf^*L^\bullet_S\otimes \omega_\pi[1]\rightarrow \pi^*L^\bullet_{\overline{\mathcal{M}}_{g}(S,\beta)/{\mathfrak{M}}_{g}}\otimes \omega_\pi[1] = \pi^{!}L^\bullet_{\overline{\mathcal{M}}_{g}(S,\beta)/{\mathfrak{M}}_{g}}.
\end{equation}  
By duality, we get a morphism
\begin{equation}\label{eq:obsthy}
    \left(R\pi_*Lf^*T^\bullet_S\right)^\vee\simeq R\pi_*\left(Lf^*L^\bullet_S\otimes \omega_\pi[1]\right)\rightarrow L^\bullet_{\overline{\mathcal{M}}_{g}(S,\beta)/{\mathfrak{M}}_{g}}.
\end{equation}

The proof that (\ref{eq:obsthy}) defines an obstruction theory is a straightforward adaptation of the proof of \cite[Proposition~6.3]{BF}.  
Note that we are not claiming that this obstruction theory is perfect.

We spell out the consequences.
Let $T$ be a $\mathbb{C}$-scheme and 
consider any family of genus $g$ stable maps
\begin{equation}\label{eq:sfamily}
    \begin{array}{ccc}
    \mathcal{C}&\stackrel{f}{\rightarrow}&S\\
    \phantom{\pi}\downarrow\pi&&\\
    T&&
    \end{array},
\end{equation}
with classifying map $g:T\to \overline{\mathcal{M}}_{g}(S,\beta)$.  Put $h=\rho\circ g:T\to{\mathfrak{M}}_{g}$.

 Let $t_0\in T$ be a closed point and let $T\subset T'$. Let the ideals of $t_0$ in $T$, 
 $t_0$ in $T'$ and $T\subset T'$ be denoted by 
 $\mathfrak{m}$ and $\mathcal{J}$ respectively.  Assume that $\mathfrak{m}\cdot \mathcal{J}=0$.
 
 Let $h':T'\to {\mathfrak{M}}_{g}$ be an extension of $h$ to $T'$ and put $C_{t_0}=\pi^{-1}(t_0)$.  Then there is an obstruction class 
\begin{equation}\label{eq:sobst}
    \alpha\in \mathbb{H}^1\left(C_{t_0},f^*T_S^\bullet\right)\otimes \mathcal{J}
\end{equation}
whose vanishing is necessary and sufficient for the family of stable maps (\ref{eq:family}) to extend to a family over stable maps over $T'$
\begin{equation}
    \begin{array}{ccc}
    \mathcal{C}'&\stackrel{f'}{\rightarrow}&S\\
    \phantom{\pi'}\downarrow\pi'&&\\
    T'&&
    \end{array},
\end{equation}
whose classifying map $g':T'\to \overline{\mathcal{M}}_{g}(S,\beta)$ satisfies $\rho\circ g' =h'$.
Furthermore, if $\alpha=0$, then the set of extensions is naturally a principal homogeneous space for $\mathbb{H}^0\left(C_{t_0},f^*T_S^\bullet\right)\otimes \mathcal{J}$.  In the above, $\mathbb{H}^p$ denotes the $p$th hypercohomology group.

However, this obstruction theory need not be perfect for general $S$.  This is the essential reason why Gromov--Witten invariants are only defined for smooth varieties.  Nevertheless, assuming that $S$ is a hypersurface in a smooth $X$, we will see that we do get an obstruction theory which is perfect, although not necessarily of amplitude $[0,1]$.

We have the triangle
\begin{equation}
    f^*(\mathcal{I}/\mathcal{I}^2)\rightarrow f^*\Omega^1_X\rightarrow Lf^*\Omega^1_S\stackrel{+}{\rightarrow},
\end{equation}
from which we deduce the triangle
\begin{equation}\label{eq:tantriangle}
  Lf^*T_S^\bullet\rightarrow f^*T_X\rightarrow f^*  \omega_S\stackrel{+}{\rightarrow}.
  \end{equation}
and finally the triangle
\begin{equation}\label{eq:compobsexplicit}
  R\pi_*Lf^*T_S^\bullet\rightarrow R\pi_*f^*T_X\rightarrow R\pi_*\left(f^*  \omega_S\right)\stackrel{+}{\rightarrow}
\end{equation}
relating the obstruction theories on the  moduli stacks of stable maps to $S$ and $X$ respectively.

Since $T_X$ and $\omega_S$ are locally free, it follows that the second and third objects in the triangle (\ref{eq:compobsexplicit}) are perfect of amplitude $[0,1]$ \cite{BehrendGW}. It then follows that first object, which describes the relative obstruction theory for $\overline{\mathcal{M}}_{g}(S,\beta)\to{\mathfrak{M}}_{g}$ is also perfect, although not necessarily of amplitude $[0,1]$.

We now use the hypothesis. We then take the (derived) restriction of (\ref{eq:tantriangle}) to $C_{t_0}$ followed by hypercohomology.  Since 
\[
H^0(C_{t_0},f^*\omega_S)=0,
\]

we get the isomorphism
\begin{equation}\label{eq:torsor}
    \mathbb{H}^0\left(C_{t_0},f^*T_S^\bullet\right)\simeq H^0(C_{t_0},f^*T_X)
\end{equation}
and the injection
\begin{equation}\label{eq:incob}
     \mathbb{H}^1\left(C_{t_0},f^*T_S^\bullet\right)\hookrightarrow H^1(C_{t_0},f^*T_X),
\end{equation}
where $\mathbb{H}^p$ denotes the $p$th cohomology object.

Now suppose that we have a family of stable maps to $X$ parametrized by a connected Noetherian scheme $B$.  

Let $T\subset B$ be the closed subscheme parametrizing those stable maps which factor through $S$, and suppose that $T$ is nonempty.  We have to show $T=B$.  Let $t_0\in T$ be a closed point with ideal $\mathfrak{m}$ and let $\mathcal{J}$ be the ideal of $T$ in $B$.  Let $T'\subset B$ be the closed subscheme defined by the ideal $\mathfrak{m}\cdot \mathcal{J}$, so that $T\subset T'$ and $T\ne T'$ if $T\ne B$ by Nakayama's Lemma.  

We show that the family of stable maps restricted to $T'$ also factors through $S$, which contradicts the definition of $T$.

From the family over $B$ restricted to $T'$ we deduce a map $h':T\to\overline{\mathcal{M}}_{g}(X,i_*\beta)$, which factors through $\overline{\mathcal{M}}_{g}(S,\beta)$. So we get obstruction classes 
\begin{equation}
    \alpha_S\in \mathbb{H}^1\left(C_{t_0},f^*T_S^\bullet\right)\otimes \mathcal{J},\ \alpha_X\in H^1\left(C_{t_0},f^*T_X\right)\otimes \mathcal{J} 
\end{equation}
with $\alpha_S$ mapping to $\alpha_X$ via (\ref{eq:incob}).  Furthermore, $\alpha_X=0$ since the original family of stable maps to $X$ over $B$ restricts to a family of stable maps to $X$ over $T'$.  We conclude 
$\alpha_S=0$ by the injectivity of (\ref{eq:incob}).  By (\ref{eq:torsor}), the set of extensions of the family over $T$ to a family of stable maps to $X$ over $'$ agrees with the set of extensions of the family over $T$ to a family of stable maps to $S$ over T'.  We conclude that the restricted family of stable maps over $T'$ factors through $S$, a contradiction.
\end{proof}

We now come to our first main result.  

\begin{thm}
    Suppose $S$ is a pre-shrinkable surface with $J=\sum_i [S_i]$.  Let $\beta\in H_2(S,\mathbb{Z})$ be any class which cannot be written as a sum of effective curve classes $\beta_i$ with $J\cdot\beta_i=\omega_S\cdot\beta_i=0$.  Then the local Gromov-Witten invariants $N^g_\beta$ of $S$ agree with the contribution of $S$ to the Gromov-Witten invariants of $X$, for any $X$ containing $S$.
\label{thm:main}
\end{thm}

In particular, we have local GW invariants for snc surfaces with ample $\omega_S^\vee$.  These surfaces are all pre-shrinkable\footnote{And are in fact shrinkable \cite{Nam}.}.

\begin{proof}
    Let $f:C\to S$ be a stable map with $f_*(C)=\beta$.  By Proposition~\ref{prop:components} we only have to show that $H^0(f^*\omega_S)=0$.  The restriction of $f^*\omega_S$ to any component $C_i$ of $C$ has nonpositive degree by hypothesis. Furthermore, we have $H^0(C_i,(f^*\omega_S)\vert_{C_i})=0$ unless $(f^*\omega_S)\vert_{C_i}$ is trivial, in which case every nonzero section of $H^0(C_i,(f^*\omega_S)\vert_{C_i})$ is nowhere vanishing.
    
 Now let $s\in H^0(f^*\omega_S)$ and let $C_i$ be a component of $C$ such that $(f^*\omega_S)\vert_{C_i}$ has negative degree.  Then $s\vert_{C_i}\equiv0$.  Similarly, the restriction of $s$ to any component $C_j$ is either identically zero or nowhere vanishing.  Since $C$ is connected, $s$ must vanish identically.
\end{proof}

We now adapt Theorem~\ref{thm:main} to the general shrinkable case.  Suppose that $S$ is shrinkable and $S\subset X$. We choose a divisor class $J=\sum_ia_i[S_i]$ that satisfies the conditions in Definition~\ref{defn:presh} with $a_i\in\mathbb{Z}$.
Let $\mathcal{I}_i$ be the ideal sheaf of $S_i$ in $X$, and let $S_J\subset X$ be the surface defined by the ideal sheaf $\prod_i\mathcal{I}_i^{a_i}$ (so that $S_J$ is non-reduced unless all $a_i=1$).  For any $\beta\in H_2(S_J,\mathbb{Z})=H_2(S,\mathbb{Z})$, we have a moduli stack $\overline{\mathcal{M}}_g(S_J,\beta)\subset \overline{\mathcal{M}}_g(X,i_*\beta)$.
Since $S_J$ is a hypersurface in $X$, Lemma~\ref{lem:perfect} shows that $R\pi_*Lf^*T_{S_J}^\bullet$ is perfect. 

\begin{ex}
    There are some pre-shrinkable surfaces for which the choice of $J$ with all $a_i=1$ does not satify (i) in Definition~\ref{defn:presh}. For instance, $J=S_1+S_2$ does not satisfy condition (i) for $S=\mathbb{F}_1\cup\mathbb{F}_7$, with a section of $\mathbb{F}_1$ of self-intersection 5 glued to the $(-7)$ curve of $\mathbb{F}_7$. Similarly, $J=S_1+S_2$ does not satisfy condition (i) for $S=\mathbb{F}_0\cup\mathbb{F}_8$ with a section of self-intersection 6 of one of the projections of $\mathbb{F}_0$ glued to the $(-8)$ curve of $\mathbb{F}_8$. One can directly check for $\mathbb{F}_1\cup\mathbb{F}_7$, $J=2\mathbb{F}_1+\mathbb{F}_7$ works and for $\mathbb{F}_0\cup\mathbb{F}_8$, $J=2\mathbb{F}_0+\mathbb{F}_8$ works. 
\label{ex:preshrink}
\end{ex}

\begin{thm}
Let $S$ be a shrinkable surface and choose a divisor class $J$ that satisfies the conditions in Definition~\ref{defn:presh}. Suppose that $\beta$ cannot be written as a sum of curve classes $[C]$ with $J\cdot C=0$.  Then $\overline{\mathcal{M}}_g(S_J,\beta)$ is a union of connected components of $\overline{\mathcal{M}}_g(X,i_*\beta)$ (as stacks), and the contribution of $S_J$ to $[\overline{\mathcal{M}}_g(X,i_*\beta)]^{\mathrm{vir}}$ is
\begin{equation}
\left(c\left(-\mathrm{ind}\left(R\pi_*Lf^*T_{S_J}^\bullet\right)-\mathrm{ind}\left(R\pi_*f^*\omega_{S_J}\right)+\mathrm{ind}\left(R\pi_*T_\pi^\bullet\right)\right)\cap
    c_F\left(\overline{\mathcal{M}}_{g}(S_J,\beta)\right)\right)_0.
 \label{eqn:ngbj}
 \end{equation}
\label{thm:generalcase}
\end{thm}

\begin{proof}
The proofs of Propositions~\ref{prop:localSiebert} and~\ref{prop:components}, and the proof of Theorem~\ref{thm:main} adapt to this situation, once we observe that $\omega_{S_J}\simeq\mathcal{O}_S(J)$.
\end{proof}

 As pointed out in Remark~\ref{rmk:shdef}, this theorem can be applied to surfaces which are not pre-shrinkable as long as they satisfy condition (i) of Definition~\ref{def:shrinkable}. For instance, $S=\mathbb{F}_0\cup\mathbb{F}_{10}$ with a section of self-intersection 8 of one of the projections of $\mathbb{F}_0$ glued to the $(-10)$ curve of $\mathbb{F}_{10}$ is not pre-shrinkable. The only choice (up to scalar multiplication) of $J$ that satisfies (i) of Definition~\ref{def:shrinkable} is $J=2\mathbb{F}_0+\mathbb{F}_{10}.$ However, $J^2\mathbb{F}_0=J^2\mathbb{F}_{10}=0$ violating (iii). Nevertheless, $S$ has local Gromov-Witten invariants and the associated physical theory has BPS invariants even though it is not a 5d SCFT.

\begin{rmk} 
This theorem suggests that unless all $a_i$ can be taken to be 1, stable maps to $S$ do not suffice for computing local Gromov-Witten invariants.  Conceivably, the schemes $S_J$ can depend on the choice of $X$ containing $S$, raising the possibility that the $N^g_\beta$ defined as the degree of (\ref{eqn:ngbj}) can depend on $X$.  
\end{rmk}

\end{document}